\documentstyle[12pt]{article}

\makeatletter
\newdimen\normalarrayskip
\newdimen\minarrayskip
\normalarrayskip\baselineskip
\minarrayskip\jot
\newif\ifold \oldtrue 
\def\arraymode{\ifold\relax\else\displaystyle\fi}
\def\eqnumphantom{\phantom{(\theequation)}}
\def\@arrayskip{\ifold\baselineskip\z@\lineskip\z@
\else
\baselineskip\minarrayskip\lineskip2\minarrayskip\fi}
\def\@arrayclassz{\ifcase \@lastchclass \@acolampacol \or
\@ampacol \or \or \or \@addamp \or
\@acolampacol \or \@firstampfalse \@acol \fi
\edef\@preamble{\@preamble
\ifcase \@chnum
\hfil$\relax\arraymode\@sharp$\hfil
\or $\relax\arraymode\@sharp$\hfil
\or \hfil$\relax\arraymode\@sharp$\fi}}
\def\@array[#1]#2{\setbox\@arstrutbox=\hbox{\vrule
height\arraystretch \ht\strutbox
depth\arraystretch \dp\strutbox
width\z@}\@mkpream{#2}\edef\@preamble{\halign
\noexpand\@halignto
\bgroup \tabskip\z@ \@arstrut \@preamble \tabskip\z@ \cr}%
\let\@startpbox\@@startpbox \let\@endpbox\@@endpbox
\if #1t\vtop \else \if#1b\vbox \else \vcenter \fi\fi
\bgroup \let\par\relax
\let\@sharp##\let\protect\relax
\@arrayskip\@preamble}
\def\eqnarray{\stepcounter{equation}%
\let\@currentlabel=\theequation
\global\@eqnswtrue
\global\@eqcnt\z@
\tabskip\@centering
\let\\=\@eqncr
$$%
\halign to \displaywidth\bgroup
\eqnumphantom\@eqnsel\hskip\@centering
$\displaystyle \tabskip\z@ {##}$%
\global\@eqcnt\@ne \hskip 2\arraycolsep
$\displaystyle\arraymode{##}$\hfil
\global\@eqcnt\tw@ \hskip 2\arraycolsep
$\displaystyle\tabskip\z@{##}$\hfil
\tabskip\@centering
&{##}\tabskip\z@\cr}
\begingroup\ifx\undefined\newsymbol \else\def\input#1 {\endgroup}\fi

\def\nl#1#2{\mathop{#2}\limits_{#1}}

\newcommand{\pl}{\partial}

\newcommand{\bqq}{\begin{equation} \label}
\newcommand{\eeq}{\end{equation}}

\begin{document}

\large
\centerline{\Large{\bf{On a generic defining function of
projective motion}}}\smallskip
\centerline{\Large{\bf{in the rigid 6-dimensional $h$-spaces
\footnote{The research was partially supported by the RFBR
grant 04-02-16538-a.}}}}

\bigskip

\centerline{{\large
Zolfira Zakirova}\footnote{Kazan State University,
e-mail: Zolfira.Zakirova@soros.ksu.ru}}

\bigskip

\abstract{\small  In this note we find a generic defining function of
projective motion in the 6-dimensional rigid $h$-space.
}

\begin{center}
\rule{5cm}{1pt}
\end{center}

\bigskip
\setcounter{footnote}{0}
\section{Introduction}

In the series of works \cite{zak0}--\cite{zak6},
we investigate the $6$-dimensional pseudo-Riemannian space $V^6(g_{ij})$
with signature $[++----]$, which admit
projective motions, i.e. continuous transformation groups preserving
geodesics. The general method of determining pseudo-Riemannian manifolds that
admit some nonhomothetic projective group $G_r$ has been developed
by A.V.Aminova \cite{am2}. A.V.Aminova has classified
all the Lorentzian manifolds
of dimension $\geq 3$ that admit nonhomothetic projective or affine
infinitesimal transformations. In each case, there were determined the
corresponding maximal projective and affine Lie algebras.
This problem is not solved for pseudo-Riemannian spaces with arbitrary signature.
In the papers \cite{zak0}--\cite{zak5}, we found 6-dimensional rigid
$h$-spaces of the $[2211]$, $[321]$, $[33]$, $[411]$ and $[51]$ types:

\noindent
{\bf {the metric of the $h$-space of the $[2211]$ type}}

\bqq{2.28}
g_{ij}dx^idx^j=e_2 (f_4-f_2)^2  \Pi_{\sigma} (f_{\sigma}-f_2) \lbrace 2A  dx^1dx^2-A^2 \Sigma_1 (dx^2)^2 \rbrace+
\eeq
$$
e_4 (f_2-f_4)^2  \Pi_{\sigma} (f_{\sigma}-f_4) \lbrace 2 \tilde{A}  dx^3dx^4-\tilde{A}^2 \Sigma_2 (dx^4)^2 \rbrace+
\sum_{\sigma} e_{\sigma} {\Pi'_i}(f_i-f_{\sigma}) (dx^{\sigma})^2,
$$
where
\bqq{2.28.4}
A=\epsilon x^1+\theta(x^2),
\quad
\tilde A=\tilde{\epsilon} x^3+\omega(x^4),
\eeq
$$
\Sigma_1=2 (f_4-f_2)^{-1}+ \sum_{\sigma}(f_{\sigma}-f_2)^{-1},
\quad
\Sigma_2=2 (f_2-f_4)^{-1}+ \sum_{\sigma}(f_{\sigma}-f_4)^{-1},
$$
$f_1=f_2=\epsilon x^2$, $f_3=f_4=\tilde{\epsilon}x^4+a$,
$\epsilon, \tilde{\epsilon}=0, 1$, $\epsilon\ne\tilde\epsilon$,
$a$ is a constant nonzero whenever $\tilde{\epsilon}=0$,
$f_{\sigma}(x^{\sigma})$, $\theta(x^2)$,  $\omega(x^4)$ are arbitrary functions, $e_i= \pm
1$.

\noindent
{\bf {the metric of the $h$-space of the $[321]$ type}}

\bqq{2.29-4}
g_{ij}dx^idx^j=e_3 (f_5-\epsilon x^3)^2(f_6-\epsilon x^3)
\lbrace (dx^2)^2 +4A  dx^1dx^3+
\eeq
$$
2(\epsilon x^1-2 A \Sigma_1)dx^2 dx^3+
((\epsilon x^1)^2-4\epsilon x^1 A \Sigma_1+
4A^2 \Sigma_3) (dx^3)^2 \rbrace+
$$
$$
e_5 (\epsilon x^3-f_5)^3
(f_6-f_5)
\lbrace 2\tilde{A}  dx^4dx^5-
\Sigma_4 {\tilde{A}}^2(dx^5)^2 \rbrace+
e_6(f_5-f_6)^2(f_5-f_6)^3
 (dx^6)^2,
$$
where
\bqq{2.29-5}
A=\epsilon x^2+\theta(x^3),
\quad
\tilde{A}=\tilde{\epsilon} x^4+\omega(x^5),
\eeq
$$
\Sigma_1=(f_6-\epsilon x^3)^{-1}+
2(f_5-\epsilon x^3)^{-1},
\quad
\Sigma_2=(f_6-\epsilon x^3)^{-2}+
2(f_5-\epsilon x^3)^{-2},
$$
$$
2\Sigma_3=(\Sigma_1)^2-\Sigma_2,
\quad
\Sigma_4=3(\epsilon x^3-f_5)^{-1}+
(f_6-f_5)^{-1},
$$
$f_1=f_2=f_3=\epsilon x^3$,
$f_4=f_5=\tilde{\epsilon}x^5+a$, $\theta(x^3)$, $\omega(x^5)$, $f_6(x^6)$ are arbitrary function,
$\epsilon, \tilde{\epsilon}=0, 1$,
$\epsilon\ne\tilde\epsilon$, $a$ is a constant nonzero whenever $\tilde{\epsilon}=0$,
$e_3, e_5, e_6 =\pm 1$.

\noindent
{\bf {the metric of the $h$-space of the $[33]$ type}}

\bqq{2.31}
g_{ij}dx^idx^j=e_3 (f_6-f_3)^3 \lbrace (dx^2)^2 +4A  dx^1dx^3+
\eeq
$$
2(\epsilon x^1-2 A \Sigma_1)dx^2 dx^3+((\epsilon x^1)^2-4\epsilon x^1
A \Sigma_1+4A^2 \Sigma_2) (dx^3)^2 \rbrace+
$$
$$
e_6 (f_3-f_6)^3
\lbrace (dx^5)^2+4 \tilde{A}  dx^4dx^6+
2(\tilde{\epsilon} x^4+2 \tilde{A} {\Sigma}_1)dx^5dx^6+
((\tilde{\epsilon} x^4)^2+4\tilde{\epsilon} x^4 \tilde{A} {\Sigma}_1+
4 \tilde{A}^2{\Sigma}_2
( (dx^6)^2 \rbrace,
$$
where
\bqq{2.35}
A=\epsilon x^2+\theta(x^3),
\quad
\tilde{A}=\tilde{\epsilon} x^4+\omega(x^6),
\eeq
$$
\Sigma_1=3 (f_6-f_3)^{-1},
\quad
\Sigma_2=3 (f_6-f_3)^{-2},
$$
$f_1=f_2=f_3=\epsilon x^3$, $f_4=f_5=f_6=\tilde{\epsilon}x^6+a$,
$\epsilon, \tilde{\epsilon}=0, 1$, $\epsilon\ne\tilde\epsilon$,
$a$ is a constant nonzero whenever $\tilde{\epsilon}=0$,
$e_i =\pm 1$,
$\theta(x^3)$, $\omega(x^6)$ are arbitrary functions.

\noindent
{\bf {the metric of the $h$-space of the $[411]$ type}}
\bqq{2.45}
g_{ij}dx^idx^j=e_4{\Pi}_{\sigma}(f_\sigma-\epsilon x^4)
\lbrace 6 Adx^1dx^4+2 dx^2dx^3+2(2\epsilon x^2-3A\Sigma_1)dx^2dx^4-
\eeq
$$
\Sigma_1(dx^3)^2+
2(\epsilon x^1-2\epsilon x^2\Sigma_1)dx^3dx^4+
4((\epsilon  x^2)^2   \Sigma_1+{\epsilon}^2
x^1x^2-\frac{3}{2}\epsilon x^1 A \Sigma_1)(dx^4)^2\rbrace+
$$
$$
3Adx^3dx^4+
12\epsilon x^2  A(dx^4)^2+\sum_{\sigma}e_{\sigma} {\Pi'}_{i}(f_i-f_\sigma)(dx^{\sigma})^2,
$$
where
\bqq{2.46}
A=\epsilon x^3+\theta(x^4),
\quad
\Sigma_1=(f_5-\epsilon    x^4)^{-1}+(f_6-\epsilon   x^4)^{-1},
\eeq
$f_1=f_2=f_3=f_4=\epsilon x^4$, $f_5(x^5)$, $f_6(x^6)$, $\theta(x^4)$ are arbitrary functions, $e_4, e_5, e_6=\pm 1$,
$\epsilon=0,1$

\noindent
{\bf {the metric of the $h$-space of the $[51]$ type}}

\bqq{2.53}
g_{ij}dx^idx^j=e_5(f_6-\epsilon x^5)
\lbrace 8 Adx^1dx^5+2 dx^2dx^4+
\eeq
$$
2(3\epsilon x^3-4A\Sigma_1)dx^2dx^5 +
(dx^3)^2-2\Sigma_1dx^3dx^4+2(2\epsilon x^2-3\epsilon x^3\Sigma_1)dx^3dx^5+
$$
$$
2(\epsilon x^1-2\epsilon x^2 \Sigma_1)dx^4dx^5+
4(3/2 \epsilon x^1\epsilon x^3+(\epsilon  x^2)^2-
2\epsilon x^1 A \Sigma_1-
$$
$$
3\epsilon x^2\epsilon x^3\Sigma_1)(dx^5)^2\rbrace+
e_6(\epsilon x^5-f_6)^5(dx^6)^2,
$$
where
\bqq{2.54}
A=\epsilon
x^4+\theta(x^5),
\quad
\Sigma_1=(f_6-\epsilon    x^5)^{-1},
\eeq
$\theta(x^5)$, $f_6(x^6)$ are arbitrary functions, $f_1=f_2=f_3=f_4=f_5=\epsilon x^5$,
$e_5, e_6=\pm 1$, $\epsilon=0, 1$.

\bigskip

The defining function of projective motion in the 6-dimensional rigid
$h$-spaces is
\bqq{2.55}
\varphi=\frac{1}{2} \sum_{i=1}^6 f_i.
\eeq

In the paper \cite{zak6}, we determined conditions of the constant curvature of the 6-dimensional
rigid $h$-spaces:

\noindent
for the $h$-space of the  $[2211]$ type
\bqq{4.47}
\rho_p-\rho_{\sigma p}=\rho_p-\rho_{pq}=\epsilon=\tilde \epsilon=0
\quad
(p\ne q, p, q=2, 4, \sigma=5, 6),
\eeq

\noindent
for the $h$-space of the  $[321]$ type
\bqq{4.47.1}
{f'}_6=\epsilon=\tilde\epsilon=0,
\eeq

\noindent
for the $h$-space of the  $[33]$ type
\bqq{4.47.2}
\epsilon=\tilde\epsilon=0,
\eeq

\noindent
for the $h$-space of the  $[411]$ type
\bqq{4.46}
\rho_p-\rho_{\sigma p}=\epsilon=\gamma_1=\gamma_2=0
\quad
(\sigma=5, 6, p=4),
\eeq

\noindent for the $h$-space of the  $[51]$ type
\bqq{4.47.3}
{f'}_6=\epsilon=0,
\eeq
where $\rho_p$, $\rho_{\sigma p}$, $\rho_{pq}$,
$\gamma_1$ and $\gamma_2$ are determined by formulas
\bqq{4.4}
\rho_p=-\frac{1}{4} \sum_{\sigma} \frac{(f'_{\sigma})^2}{(f_{\sigma}-f_p)^2
g_{\sigma \sigma}},
\quad
\rho_{pq}=-\frac{1}{4} \sum_{\sigma} \frac{(f'_{\sigma})^2}
{(f_{\sigma}-f_p)(f_{\sigma}-f_q)g_{\sigma \sigma}},
\eeq
$$
\rho_{\sigma p}=-\frac{1}{4}\frac{(f'_{\sigma})^2}{(f_{\sigma}-f_p)
g_{\sigma \sigma}} \lbrace \frac{2 {f''}_{\sigma}}{(f'_{\sigma})^2}-
\frac{1}{f_{\sigma}-f_p}+
\sum_{i, i \ne \sigma} (f_i-f_{\sigma})^{-1} \rbrace
$$
$$
-\frac{1}{4} \sum_{\gamma, \gamma \ne \sigma} \frac{(f'_{\gamma})^2}
{(f_{\gamma}-f_p)(f_{\gamma}-f_{\sigma})g_{\gamma \gamma}}.
$$
\bqq{4.23.1}
\gamma_1=-\frac{1}{4} \sum_{\sigma} \frac{(f'_{\sigma})^2}
{(f_{\sigma}-f_4)^3 g_{\sigma \sigma}},
\quad
\gamma_2=-\frac{1}{4} \sum_{\sigma} \frac{(f'_{\sigma})^2}
{(f_{\sigma}-f_4)^4 g_{\sigma \sigma}}.
\eeq

The results obtained in the  papers \cite{zak0}-\cite{zak6} are necessary for investigating
properties of the defining function of projective motion in the rigid 6-dimensional
$h$-spaces.

\section{Generic defining function of projective motion in the rigid 6-dimensional $h$-spaces}

In this section we prove the following theorem that manifestly gives any
defining function of projective motion.

\bigskip

\noindent
{\bf Theorem. } {\it Any defining function of
projective motion in the rigid $h$-spaces
of nonconstant curvature can be presented as $\phi=a_1 \varphi$,
where $a_1$ is a constant and
$\varphi$ is determined by {\rm (\ref{2.55})}.}

\noindent
{\bf Proof.} Suppose we are given with a vector field ${\xi}^i$
that gives the projective
motion with the defining function $\phi$ in the
rigid $h$-spaces.
Then, for the tensor $b_{ij}={\nl{\xi}L} g_{ij}$ the following
Eisenhart equations \cite{ezen1} are fulfilled
\bqq{14}
b_{ij,k}=2g_{ij} \phi_{,k}+g_{ik} \phi_{,j}+
g_{jk} \phi_{,i}.
\eeq
Their integrability conditions are
\bqq{15}
b_{mi} R_{jkl}^m+b_{mj} R_{ikl}^m=g_{ik} \phi_{,jl}+g_{jk}\phi_{,il}-g_{li}
\phi_{,jk}-g_{lj} \phi_{,ik}.
\eeq
Let us show that for any projective motion in the
rigid $h$-spaces the following conditions are satisfied
\bqq{16}
b_{\alpha \beta}=0,
\quad
\phi_{,\alpha \beta}=0,
\eeq
$\alpha$, $\beta$ being indices of vanishing components of the metric $g_{\alpha \beta}$
of the $h$-spaces of the $[2211]$, $[321]$, $[33]$, $[411]$ and $[51]$ type.

For the $h$-space of the  $[2211]$ type, the set $(ijkl)=(3112)$, $(3314)$,
$(3\sigma 1\sigma)$, $(1334)$, $(1132)$, $(1\sigma 3\sigma)$. From (\ref{15}), one finds
\bqq{17}
b_{13}R_{112}^1=-g_{12} \phi_{,13},
\quad
b_{13}R_{314}^1=-g_{34} \phi_{,13},
\quad
b_{13}R_{\sigma 1 \sigma}^1=-g_{\sigma \sigma} \phi_{,13},
\eeq
\bqq{18}
b_{13}R_{334}^3=-g_{34} \phi_{,13},
\quad
b_{13}R_{132}^3=-g_{12} \phi_{,13},
\quad
b_{13}R_{\sigma 3 \sigma}^3=-g_{\sigma \sigma} \phi_{,13},
\eeq
From (\ref{17}) it follows that
$$
b_{13}(\chi_2-\rho_{24})=b_{13}
(\chi_2-\rho_{\sigma 2})=0,
$$
and  from (\ref{18})
$$
b_{13}(\chi_4-\rho_{24})=
b_{13}(\chi_4-\rho_{\sigma 4})=0.
$$
If $b_{13} \ne 0$, then
$\chi_p-\rho_{\sigma p}=0$ $(p=2, 4)$. Since we consider $h$-spaces of the nonconstant curvature, one
necessarily puts $b_{13}=0$, hence $\phi_{,13}=0$. Similarly one can get other equalities (\ref{16}).

In the case of the $h$-space of the $[321]$ type, when $(ijkl)=(1212)$, $(1415)$,
$(1515)$, $(1616)$, $(2212)$, $(1123)$, $(2415)$, $(4125)$, equations
(\ref{15}) leads to the following identities
\bqq{19}
b_{11} R_{212}^1=-g_{22} \phi_{,11},
\quad b_{11} R_{123}^1=b_{11} R_{425}^1=0,
\eeq
\bqq{20}
b_{11} \frac{R_{415}^1}{g_{45}}=
b_{11}\frac{R_{515}^{1}}{g_{55}}=
b_{11}\frac{R_{616}^1}{g_{66}}=\phi_{,11},
\eeq
hence
$$
b_{11}(\chi_3-\rho_{63})=b_{11}(\chi_3-\rho_{53})=
b_{11}(\chi_3-\rho_{53}-
\frac{\tilde \epsilon \omega'}{(f_5-f_3) \tilde A g_{55}})=0.
$$
If $b_{11} \ne 0$, then $\chi_3-\rho_{63}=\chi_3-\rho_{53}=\tilde \epsilon=
0$. From  this it follows that $\epsilon=\rho_3-\rho_{53}=0$, and $f'_6=0$. Therefore,
$\chi_5-\rho_{53}=\chi_5-\rho_{65}=\gamma=0$, this leads to the constant curvature $h$-space. In such a manner,
$b_{11}=\phi_{,11}=0$. Similarly one can get other equalities (\ref{16}).

To prove that
$b_{11}=\phi_{,11}=0$ for the $h$-space of the $[33]$ type, set $(ijkl)=(1212)$, $(1416)$ and $(1516)$
in (\ref{15}). Then, one obtains
\bqq{21}
b_{11} R_{212}^1=-g_{22} \phi_{,11},
\quad
b_{11} R_{416}^1=-g_{46} \phi_{,11},
\eeq
\bqq{22}
b_{11}(\frac{R_{212}^1}{g_{22}}-\frac{R_{516}^1}{g_{56}})=0.
\eeq
Since $R_{416}^1=0$ and $g_{46} \ne 0$, we get $\phi_{,11}=0$.
Then, from (\ref{22}) it follows that $b_{11}=0$, since
$R_{212}^1 \ne 0$.
Similarly one can get other equalities (\ref{16}).

Putting $(ijkl)=(1213),(1313), (1\sigma 1 \sigma), (3312), (1314)$ in the
case of the $h$-space of the $[411]$ type, one gets
\bqq{23}
b_{11} R_{213}^1=-g_{23} \phi_{,11},
\quad
b_{11} R_{313}^1=-g_{33} \phi_{,11},
\eeq
\bqq{24}
b_{11} R_{\sigma 1 \sigma}^1=-g_{\sigma \sigma} \phi_{,11},
\eeq
\bqq{25}
b_{13} R_{312}^1=-g_{23} \phi_{,13},
\quad
b_{11} R_{314}^1+b_{13} R_{114}^1=-g_{14} \phi_{,13}-g_{34} \phi_{,11}.
\eeq
Since
$$
R_{213}^1=R_{312}^1=g_{23} \rho_4,
\quad
R_{114}^1=g_{14} \rho_4,
$$
$$
R_{313}^1=g_{33}\rho_4+\gamma_1 g_{23}+\frac{2 {\epsilon}^2}{9 A^2},
\quad
R_{\sigma 1 \sigma}^1=g_{\sigma \sigma }\rho_{\sigma 4},
$$
$$
R_{314}^1=g_{34} \rho_4+\gamma_1 g_{24}+
\gamma_2 g_{14}+\frac{2\epsilon}{3 A^2}(\theta'-\frac{4{\epsilon}^2 x^2}
{3}),
$$
then, from (\ref{23}) and (\ref{24}), we have
$$
b_{11}(\rho_4-\rho_{\sigma 4})=b_{11}(\rho_4-\rho_{\sigma 4}+
\gamma_3 \frac{g_{23}}{g_{33}}+\frac{2 {\epsilon}^2}{9 A^2 g_{23}})=0.
$$
From (\ref{24}) and (\ref{25}), one finds
$$
b_{11}(\rho_4-\rho_{\sigma 4}+
\gamma_3 \frac{g_{24}}{g_{34}}+\gamma_4 \frac{g_{14}}{g_{34}}+
\frac{2 \epsilon}{3 A^2 g_{34}}(\theta'-\frac{4 {\epsilon}^2 x^2}{3}))=0.
$$
Let us $b_{11} \ne 0$. From the last relations it follows that $\rho_4-\rho_{\sigma 4}=\epsilon=\gamma_1=\gamma_2=0$,
which are necessary and sufficient
conditions of the constant curvature of these $h$-spaces, therefore, $b_{11}=0$,
and, hence, $\phi_{,11}=0$. Similarly one can get other equalities (\ref{16}).

For the $h$-space of the $[51]$ type from  integrability conditions
one can similarly get equalities (\ref{16}).

Consider equation (\ref{14}). Putting
$(ijk)=(1\sigma\sigma), (3\sigma\sigma), (\sigma\tau\sigma),
(\tau\sigma\tau)$ $(\sigma, \tau=5, 6, \sigma \ne \tau)$ for the $h$-space of the $[2211]$ type, we determine
\bqq{26}
\phi_{,1}=0,
\quad
\phi_{,3}=0,
\eeq
\bqq{27}
\phi_{,\sigma}=f'_{\sigma} P_{\sigma\tau},
\quad
\phi_{,\tau}=f'_{\tau}P_{\sigma\tau},
\quad
P_{\sigma\tau}\equiv \frac{1}{2}
\frac{b_{\tau\tau}g_{\tau\tau}^{-1}-b_{\sigma\sigma}g_{\sigma\sigma}^{-1}}{f_\tau-f_\sigma}, (\sigma\ne\tau).
\eeq
From this it follows that $f'_{\sigma} \phi_{,\tau}=f'_{\tau} \phi_{,\sigma}$.
Therefore, from the condition  $\phi_{,\sigma\tau}=0$, one gets
\bqq{28}
\pl_{\sigma\tau}\phi=0.
\eeq
Using relations obtained from (\ref{14}) with
$(ijk)= (1\tau2), (2\tau\tau), (\tau\tau\tau)$, we
find
\bqq{29}
\phi_{,\tau}=\frac{1}{2} \frac{f'_{\tau}}{f_\tau-f_2}(\frac{b_{\tau \tau}}{g_{\tau\tau}}-\frac{b_{12}}{g_{12}}),
\eeq
\bqq{30}
\phi_{,2}=\frac{\epsilon}{f_\tau-f_2}(\frac{b_{\tau \tau}}{g_{\tau\tau}}-\frac{b_{12}}{g_{12}}),
\eeq
\bqq{31}
\frac{\pl_{\tau}b_{\tau \tau}}{g_{\tau\tau}}=-f'_\tau
{\nl{i, j \ne \tau}\sum}(f_i-f_\tau)^{-1}\frac{b_{\tau\tau}}{g_{\tau\tau}}+4\phi_{,\tau}.
\eeq
From (\ref{29}) and (\ref{30}), we have
\bqq{47.1}
2\epsilon \phi_{,\tau}=f'_{\tau} \phi_{,2}.
\eeq
These relations and conditions $\phi_{,2 \tau}=0$ give
\bqq{47.2}
\pl_{2\tau} \phi=0.
\eeq
Differentiating (\ref{29}) w.r.t. $x^{\tau}$ and taking into account
(\ref{31}), one obtains
\bqq{47.3}
f'_{\tau} \pl_{\tau}\phi_{,\tau}=f''_{\tau} \phi_{,\tau}.
\eeq
Now, from (\ref{14}) with $(ijk)=(142), (234)$, one finds
\bqq{47.7}
\phi_{,2}=\frac{\tilde\epsilon}{f_4-f_2}
(\frac{b_{34}}{g_{34}}-\frac{b_{12}}{g_{12}}),
\quad
\phi_{,4}=\frac{\epsilon}{f_4-f_2}
(\frac{b_{34}}{g_{34}}-\frac{b_{12}}{g_{12}}),
\eeq
hence, $\epsilon \phi_{,4}=\tilde \epsilon \phi_{,2}$ and, since
$\phi_{,24}=0$, one obtains $\pl_{24} \phi=0$.

Assuming $(ijk)=(3\tau 4), (4\tau\tau)$  in equation (\ref{14}), one obtains
\bqq{47.8}
2\tilde \epsilon \phi_{, \tau}=f'_{\tau} \phi_{,4},
\eeq
therefore, taking into account $\phi_{,4\tau}=0$, we have $\pl_{4\tau}\phi=0$.
If not all $f'_{\tau}=0$, then, from (\ref{47.1}),
(\ref{47.3}) and (\ref{47.8}), one finds
\bqq{47.9}
\phi_{, \tau}=\frac{1}{2} a_1 f'_{\tau},
\quad
\phi_{,2}=a_1 \epsilon,
\quad
\phi_{,4}=a_1 \tilde \epsilon.
\eeq
The first equation of (\ref{47.9}) is correct whenever $f'_{\tau}=0$ due to
relation (\ref{29}).
Integrating equation (\ref{47.9}), one obtains
\bqq{48.0}
\phi=\frac{1}{2} a_1 \nl{i=1}{\sum}^6 f_i=a_1 \varphi,
\quad
f_1=f_2=\epsilon x^2,
\quad
f_3=f_4=\tilde\epsilon x^4+a,
\eeq
where $\varphi$ determined by formula (\ref{2.55}).
If all $f'_{\tau}=0$, then, from (\ref{15}) with $(ijkl)=
(\tau 2\tau 2), (\tau 4 \tau 4)$, one finds
$\pl_{22}\phi=\pl_{44}\phi=0$, and, therefore, (\ref{48.0}) is correct.

For the $h$-space of the  $[321]$ type, putting in equation (\ref{14})
$(ijk)=(166), (266), (466)$, one obtains $\phi_{,1}=\phi_{,2}=\phi_{,4}=0$.

Equation (\ref{14}) with $(ijk)=(163), (366), (666)$ leads to
\bqq{48.01}
\phi_{,6}=\frac{1}{2}\frac{f'_{6}}{f_6-f_3}
(\frac{b_{66}}{g_{66}}-\frac{b_{13}}{g_{13}}),
\eeq
\bqq{48.02}
\phi_{,3}=\frac{3}{2}\frac{\epsilon}{f_6-f_3}
(\frac{b_{66}}{g_{66}}-\frac{b_{13}}{g_{13}}),
\eeq
\bqq{48.03}
\frac{\pl_{6}b_{66}}{g_{66}}=-f'_6
{\nl{i, j \ne 6}\sum}(f_i-f_6)^{-1}\frac{b_{66}}{g_{66}}+
4\phi_{,6},
\eeq
hence
\bqq{48.1}
f'_{6} \pl_{6} \phi_{,6}=f''_{6} \phi_{,6},
\eeq
\bqq{48.2}
3 \epsilon \phi_{,6}=f'_{6} \phi_{,3}.
\eeq
Since $\phi_{,36}=0$, then, from equation (\ref{48.2}), it follows that
$\pl_{36}\phi=0$.
In case of $(ijk)=(465), (566)$, from (\ref{14}) we have
\bqq{48.21}
\phi_{,6}=\frac{1}{2} \frac{f'_{6}}{f_6-f_2}(\frac{b_{66
}}{g_{66}}-\frac{b_{12}}{g_{12}}),
\eeq
\bqq{48.22}
\phi_{,2}=\frac{\epsilon}{f_6-f_2}(\frac{b_{66
}}{g_{66}}-\frac{b_{12}}{g_{12}}),
\eeq
and, therefore,
\bqq{48.23}
2\epsilon \phi_{,6}=f'_{6} \phi_{,2}.
\eeq
Substituting $(ijk)=(345), (513)$ in (\ref{14}), one finds
\bqq{48.3}
\phi_{,3}=\frac{3}{2}\frac{\epsilon}{f_5-f_3}
(\frac{b_{45}}{g_{45}}-\frac{b_{13}}{g_{13}}),
\quad
\phi_{,5}=\frac{\tilde\epsilon}{f_5-f_3}
(\frac{b_{45}}{g_{45}}-\frac{b_{13}}{g_{13}}),
\eeq
therefore,
\bqq{48.4}
3\epsilon \phi_{,5}=2\tilde\epsilon \phi_{,3}.
\eeq
From this relation, taking into account $\phi_{,35}=0$, we find
$\pl_{35} \phi=0$. If $f'_6 \ne 0$, then, integrating the equations
(\ref{48.1}), (\ref{48.2}) and (\ref{48.4}), one obtains
\bqq{48.211}
\phi=\frac{1}{2} a_1 \nl{i=1}{\sum}^6 f_i=a_1 \varphi,
\quad
f_1=f_2=f_3=\epsilon x^3,
\quad
f_4=f_5=\tilde\epsilon x^5+a,
\eeq
where $\varphi$ is determined by formula (\ref{2.55}).

On the contrary, from equation (\ref{15}) with $(ijkl)=(3312), (3636)$,
$(6565)$ one can obtain $\pl_{33} \phi=\pl_{55} \phi=0$ which leads to the
desired result.

From (\ref{14}) with $(ijk)=(166), (266), (433), (533), (346),
(613)$ in the case of the $h$-space of the $[33]$ type, one obtains the following relations
$$
\phi_{,1}=\phi_{,2}=\phi_{,4}=\phi_{,5}=0,
$$
$$
\phi_{,3}=\frac{3}{2}\frac{\epsilon}{f_6-f_3}
(\frac{b_{46}}{g_{46}}-\frac{b_{13}}{g_{13}}),
\quad
\phi_{,6}=\frac{3}{2}\frac{\tilde\epsilon}{f_6-f_3}
(\frac{b_{46}}{g_{46}}-\frac{b_{13}}{g_{13}}).
$$
Hence
\bqq{48.5}
\tilde\epsilon \phi_{,3}=\epsilon \phi_{,6}.
\eeq
Therefore, from the condition $\phi_{,36}=0$, one obtains $\pl_{36} \phi=0$.
If $\tilde\epsilon \ne 0$, then, differentiating relation (\ref{48.5})
w.r.t.  $x^3$, we have $\pl_{33} \phi=0$. If $\tilde\epsilon=0$, then
$\epsilon \ne 0$, since we consider the space of nonconstant curvature.
In case when  $\tilde\epsilon=0$ and $\epsilon \ne 0$,
the relation $\pl_{33}\phi=0$ follows from
(\ref{15}) with $(ijkl)=(3312), (3436)$.  Similarly one can get that
$\pl_{66} \phi=0$. Integrating the obtained equations, one finds
\bqq{48.6}
\phi=\frac{1}{2} a_1 \nl{i=1}{\sum}^6 f_i=a_1 \varphi,
\quad
f_1=f_2=f_3=\epsilon x^3,
\quad
f_3=f_4=f_6=\tilde\epsilon x^6+a,
\eeq
where
$\varphi$ is determined by formula (\ref{2.55}).

For the $h$-space of the $[411]$ type from (\ref{14}) with
$(ijk)=(1\tau\tau), (2\tau\tau),(3\tau\tau)$ $(\tau=5, 6)$, one obtains
$$
\phi_{,1}=\phi_{,2}=\phi_{,3}=0.
$$
For this $h$-space relations (\ref{27}) are also correct. From these it follows that $\pl_{\sigma \tau}\phi=0$
$(\sigma, \tau=5, 6, \sigma \ne \tau)$.  Using the obtained results, from
(\ref{14}) á $(ijk)=(1\tau4), (4\tau\tau), (\tau\tau\tau)$, we find
\bqq{48.7} \phi_{,\tau}=\frac{1}{2} \frac{f'_{\tau}}{f_\tau-f_4}
(\frac{b_{\tau \tau}}{g_{\tau\tau}}-\frac{b_{14}}{g_{14}}),
\eeq
\bqq{48.8}
\phi_{,4}=2\frac{\epsilon}{f_\tau-f_4}
(\frac{b_{\tau \tau}}{g_{\tau\tau}}-\frac{b_{14}}{g_{14}}),
\eeq
\bqq{48.9}
\frac{\pl_{\tau}b_{\tau \tau}}{g_{\tau\tau}}=-f'_\tau
{\nl{i, j \ne \tau}\sum}(f_i-f_\tau)^{-1}\frac{b_{\tau\tau}}
{g_{\tau\tau}}+4\phi_{,\tau},
\eeq
hence,
\bqq{49.0}
f'_{\tau} \pl_{\tau} \phi_{, \tau}=f''_{\tau} \phi_{, \tau},
\eeq
\bqq{49.1}
4\epsilon \phi_{, \tau}=f'_{\tau} \phi_{,4}.
\eeq
Using the condition $\phi_{,4\tau}=0$ and relation (\ref{49.1}), one gets
$\pl_{4\tau}\phi=0$. If not all $f'_{\tau} \ne 0$, then, from obtained formulas it follows that
\bqq{49.11}
\phi=\frac{1}{2} a_1 \nl{i=1}{\sum}^6 f_i=a_1 \varphi,
\quad
f_1=f_2=f_3=f_4=\epsilon x^4,
\eeq
where $\varphi$ is determined by formula (\ref{2.55}).
In the opposite case, when all $f'_{\tau}=0$, we have
$$
\phi_{,\tau}=\phi_{,14}=\phi_{,23}=\phi_{,24}=\phi_{,\tau \tau}=0,
$$
\bqq{49.2}
\phi_{,34}=-\frac{2 {\epsilon}^2}{3A(f_{\tau}-f_4)}
(\frac{b_{\tau\tau}}{g_{\tau \tau}}-\frac{b_{14}}{g_{14}}),
\eeq
$$
\phi_{,44}=\pl_{44}\phi-\frac{2 \epsilon(3 \theta'-2{\epsilon}^2 x^2)}
{3A(f_{\tau}-f_4)}
(\frac{b_{\tau\tau}}{g_{\tau \tau}}-\frac{b_{14}}{g_{14}}),
$$
where $\epsilon=1$, since when  $\epsilon=f'_{\tau}=0$ the space is of constant curvature.

Put $(ijkl)=(4413), (4\tau 4\tau)$ in equation (\ref{15}). One
finds
\bqq{49.3}
b_{14} R_{413}^1+b_{24} R_{413}^2=g_{14} \phi_{,34},
\eeq
\bqq{49.4}
b_{14}R_{\tau 4 \tau}^1+ b_{24} R_{\tau 4\tau}^2+b_{\tau\tau}
R_{44\tau}^{\tau}=-g_{\tau\tau} \phi_{,44}.
\eeq
Expressing $b_{24}$ from (\ref{49.3}) and substituting the result in (\ref{49.4}),
taking into account (\ref{49.2}), one obtains $\pl_{44} \phi=0$ coming to the result needed.

Using equations (\ref{14}) with $(ijk)=(166), (266), (366), (466), (165), (566)$,
$(666)$ in the case of the $h$-space of the $[51]$ type, one gets
$$
\phi_{,1}=\phi_{,2}=\phi_{,3}=\phi_{,4}=0,
$$
\bqq{49.5}
\phi_{,6}=\frac{1}{2} \frac{f'_6}{f_6-f_5}
(\frac{b_{66}}{g_{66}}-\frac{b_{15}}{g_{15}}),
\quad
\phi_{,5}=\frac{5}{2}\frac{\epsilon}{f_6-f_5}
(\frac{b_{66}}{g_{66}}-\frac{b_{15}}{g_{15}}),
\eeq
$$
\frac{\pl_{6}b_{66}}{g_{66}}=-5f'_6
(f_5-f_6)^{-1}\frac{b_{66}}
{g_{66}}+4\phi_{,6},
$$
and, hence,
\bqq{49.6}
f'_{6} \pl_{6} \phi_{,6}=f''_{6} \phi_{,6},
\quad
5\epsilon \phi_{,6}=f'_{6} \phi_{,5}.
\eeq
The last equation and the condition $\phi_{,56}=0$ lead to
\bqq{49.7}
\pl_{56} \phi=0.
\eeq
Integrating equations (\ref{49.6}) and (\ref{49.7}) when $f'_{6} \ne 0$,
one obtains
\bqq{49.8}
\phi=\frac{1}{2} a_1 \nl{i=1}{\sum}^6 f_i=a_1 \varphi,
\quad
f_1=f_2=f_3=f_4=f_5=\epsilon x^5,
\eeq
where $\varphi$ is obtained by formula (\ref{2.55}).

In the case of $f'_{6}=0$, from the formulas obtained we can find
$$
\phi_{,6}=\phi_{,15}=\phi_{,24}=\phi_{,25}=\phi_{,33}=\phi_{,34}=
\phi_{,35}=\phi_{,66}=0,
$$
\bqq{49.9}
\phi_{,45}=-\frac{15 {\epsilon}^2}{16A(f_{6}-f_5)}
(\frac{b_{66}}{g_{66}}-\frac{b_{15}}{g_{15}}),
\eeq
$$
\phi_{,55}=\pl_{55}\phi-\frac{5 \epsilon(4 \theta'-3{\epsilon}^2 x^3)}
{4A(f_{6}-f_5)}
(\frac{b_{66}}{g_{66}}-\frac{b_{15}}{g_{15}}),
$$
where $\epsilon=1$, since we consider the space of nonconstant curvature.

From (\ref{15}) with $(ijkl)=(5514), (5656)$ and
taking into account (\ref{14}), one obtains the equality $\pl_{55}\phi=0$, which leads to relation (\ref{49.8}).

Therefore, the theorem is proved.

\bigskip

\bigskip
I am grateful to professor A.V.Aminova for constant attention to this work
and for useful discussions.


\begin{thebibliography}{30}

\bibitem[A95]{am2}Aminova A. V. Lie algebras of infinitesimal projective
transformations of Lorentz manifolds (in Russian) // Uspekhi Mat. Nauk
{\bf 50} (1995), no. 1(301), 69-142.

\bibitem[E]{ezen1} L.P. Eizenhart, Riemannian geometry, Princeton University Press,
1997.

\bibitem[Z]{zak0} Z. Zakirova, PhD thesis, 2001, Kazan State University

\bibitem[Z1]{zak1} Z. Zakirova, First integral of the geodesics equations in
6d h-spaces // in: Proc. of Geometry Seminar, Kazan, 1996

\bibitem[Z2]{zak2} Z. Zakirova, 6-dimensional $h$-spaces of the special
type // Lobachevskiy Int. geom. seminar, Kazan, 1997:
Abstracts of contributed papers -- Kazan, 1997.- P. 52.

\bibitem[Z3]{zak3} Z. Zakirova, On one solution of the Eisenhart
equation. // "Latest problems of the field theory. 1999-2000".
Editor  A. V. Aminova.- Kazan, 2000. '. 104-109.

\bibitem[Z4]{zak4} Z. Zakirova, First integral of the geodesics equations
in  $h$-space of the [411] type. // Izv. Vuzov. Mathematics.
N9 (448), 1999, '. 78-79.

\bibitem[Z5]{zak5} Z.Zakirova, On projective group properties of the
6d h-space of the type [33], math.DG/0101052; submitted to
"Classical \& Quantum Gravity"

\bibitem[Z6]{zak6}Z.Zakirova, Rigid 6-dimensional h-spaces of constant curvature, math.DG/0301105



\end{thebibliography}
\end{document}